\documentclass[final,1p,times,sort]{elsarticle}
\usepackage{amsmath,amsfonts,amssymb,amsthm}
\usepackage[all]{xy}
\usepackage{graphicx,subfig,aliascnt,hyperref}
\usepackage{xcolor}

\journal{arXiv}

\theoremstyle{plain}
\newtheorem{theorem}{Theorem}
\newtheorem{proposition}[theorem]{Proposition}

\newtheorem{corollary}[theorem]{Corollary}

\theoremstyle{definition}
\newtheorem{definition}[theorem]{Definition}

\newtheorem{example}[theorem]{Example}
\newtheorem{remark}[theorem]{Remark}
\newtheorem{notation}[theorem]{Notation}

\DeclareMathOperator{\Hom}{Hom}
\DeclareMathOperator{\rk}{rk}
\DeclareMathOperator{\Spec}{Spec}
\DeclareMathOperator{\Cont}{Cont}

\begin{document}

\begin{frontmatter}

\title{On the representation of measures over bounded lattices}

\author[dm]{C\'esar Massri\corref{correspondencia}}
\ead{cmassri@caece.edu.ar}

\author[df]{Federico Holik\fnref{financiado}}

\address[dm]{Department of Mathematics, CAECE, Argentina.}
\cortext[correspondencia]{Address for correspondence:
IMAS Institute, UBA-CONICET, Buenos Aires, Argentina.
Postal address: C1428EGA. Phone number: 54 011 4576-3390 (int. 915).}

\address[df]{Universidad Nacional de La Plata, Instituto
de F\'{\i}sica (IFLP-CCT-CONICET), C.C. 727, 1900 La Plata, Argentina
\\
Center Leo Apostel for Interdisciplinary Studies and, Department
of Mathematics, Brussels Free University Krijgskundestraat 33, 1160 Brussels, Belgium.}

\fntext[financiado]{The authors were fully supported by CONICET.}

\begin{abstract}
In this paper we investigate measures over bounded lattices,
extending and giving a unifying treatment to previous works. In
particular, we prove that the measures of an arbitrary bounded
lattice can be represented as measures over a suitably chosen
Boolean lattice. Using techniques from algebraic geometry, we also
prove that given a bounded lattice $X$ there exists a scheme
$\mathcal{X}$ such that a measure over $X$ is the same as a
(scheme-theoretic) measure over $\mathcal{X}$. We also define the
measurability of a lattice, and describe measures over finite
lattices.
\end{abstract}

\end{frontmatter}

\section{Introduction}

\subsection*{Overview of the problem}

Different notions of measure theory over a $\sigma$-algebra have been
adapted to measures over Boolean and orthomodular lattices (see for
example \cite{MR1670882,MR1662485,MR0584017}). Many important
results have been generalized to these settings, such as extension,
convergence and decomposition theorems, etc. (see
\cite[Ch.6]{MR2015280},
\cite{MR0304608,MR0503098,MR534996,MR867884,MR1000195,MR854631,MR1018455,MR1088647,MR1111759} and
\cite{MR1670882,MR0513905}). Due to Stone's representation theorem,
\cite{MR1501865}, measures over Boolean lattices are known. In this
paper we focus on measures over bounded lattices.

The definition of measure over a lattice depends on the type of the
particular lattice. For example, if the lattice is not bounded
below, the axiom $\nu(0)=0$ is meaningless. Also, for
$\sigma$-measures we need a $\sigma$-lattice. But, if the lattice
$\mathcal{L}$ is $\sigma$-orthomodular, we have the following
definition: a function $\rho:\mathcal{L}\rightarrow\mathbb{R}$ is
called a \emph{$\sigma$-orthomodular-measure} if $\rho(0)=0$ and
$\rho(\vee x_i)=\sum \rho(x_i)$ for any countable collection of
pairwise orthogonal elements $\{x_i\}\subseteq\mathcal{L}$. From
this definition, it follows the definition of orthomodular-measures
and ($\sigma$-) Boolean-measures (Boolean lattices are
orthomodular). Orthomodular (and $\sigma$-orthomodular) measures
find applications as states over generalized physical theories
\cite{MR0274355,Holik2014293,Holik-Massri-Plastino}.

In the case of a measure over a Boolean lattice, we can replace the
definition by an equivalent one. A \emph{Boolean-measure} over a
Boolean lattice $B$ is a function $\nu:B\rightarrow \mathbb{R}$
satisfying the inclusion-exclusion principle, that is, $\nu(0)=0$
and $\nu(x\vee y)+\nu(x\wedge y)=\nu(x)+\nu(y)$. The advantage of
this definition is that we do not need to refer to the
orthocomplementation and we can extend it to distributive lattices
(which are not necessarily Boolean). A \emph{distributive-measure}
over a distributive lattice $D$ is a function $\nu:D\rightarrow
\mathbb{R}$ satisfying the inclusion-exclusion principle. Using the
distributivity law on $D$, it follows the inclusion-exclusion
principle for any number of elements,

\[
\nu(x_1\vee\dots\vee x_k)=
\sum_{s=1}^k(-1)^{s+1}\sum_{1\leq i_1<\dots<i_s\leq k}
\nu(x_{i_1}\wedge\dots\wedge x_{i_s}).
\]

Notice that the definition of distributive-measure is compatible
with the construction that embeds a distributive lattice into a
Boolean lattice by adding formal complements (see \cite{MR1564003}).
Given an arbitrary bounded lattice $X$, we call a function
$\nu:X\rightarrow\mathbb{R}$ a \emph{bounded-measure}, if it
satisfies the inclusion-exclusion principle for any finite number of
elements. As it will be clear below, the notion of bounded measure
turns out to be the same as that of distributive valuation (see for
example, \cite{geis1}, page 465). Also, notice that, when $X$ is
Boolean, the definitions of bounded-measure and Boolean-measure
coincide. However, if, for example, $X$ is orthomodular,
bounded-measures are different from orthomodular-measures. In fact,
any orthomodular-measure can be obtained by pasting Boolean-measures
\cite[p.127]{MR0274355}. Specifically, we assign the Boolean-measure
$\rho|_{B}:B\rightarrow \mathbb{R}$ to a given orthomodular-measure
$\rho:X\rightarrow\mathbb{R}$, where $B\subseteq X$ is a Boolean
sublattice. Recall that any orthomodular lattice is the union of its
maximal Boolean sublattices \cite{MR0432504}. Finally, regarding
measures over a $\sigma$-lattice, in a future paper we plan to apply
the definition of bounded-measures and density arguments, to treat
(and even characterize) $\sigma$-bounded-measures. In this paper we
concentrate in bounded measures (see \cite{PaperGeneralizedStates}
for the case of orthomodular measures). It is important to mention
that bounded measures are of interest in many branches of
mathematics. Remarkably enough, they appear as dimension functions
in the problem of classification of factor von Neumann algebras (see
Remark \ref{rem:DimensionFunctions} in Section \ref{s:Examples}).

\subsection*{Previous work}

The definition of measure dates back to the early twentieth century
in the work of Hausdorff and Carath\'eodory. One of the first
authors who extended the notion of measure to lattices is M. F.
Smiley \cite{MR0001445}. The theory has substantially grown since
then, and the influence of quantum mechanics inspired many
developments in the non-Boolean setting \cite{MR0274355}. Particular
efforts were dedicated to find results about existence of measures
over orthomodular lattices (see
\cite{MR0689178,MR0130703,MR2015280,MR0274355}). It is also worth to
mention here two salient results that characterize measure spaces.
The first one is Radon-Nikodym's theorem in measurable spaces
\cite[\S 9]{MR0466463}. The second one is Gleason's theorem in the
lattice of projectors on a separable Hilbert space and its
generalizations (see \cite{MR0096113, MR2015280} and also
\cite{bloch1946nuclear} for a representation in dimension two).

In this work we study measures over bounded lattices (i.e., bounded
measures) by elaborating on the works of G.C. Rota
\cite{rota-valuation-ring} and L. Geissinger \cite{geis1}. In
\cite{rota-valuation-ring}, Rota defined, for each distributive
lattice $X$, a ring $V(X)$ called \emph{valuation ring}, which is
defined as $\mathbb{Z}^{\oplus X}/J$, being $J$ the submodule
generated by $x\vee y +x\wedge y-x-y$ for all $x,y\in X$. There,
Rota proved that the valuation ring represents measures over
distributive lattices. In \cite{rota-measure-logical}, the authors
study measures over Boolean lattices, and prove that given a free
Boolean lattice $B=B(P)$ over a set $P$, then $V(B)$
is isomorphic to $\mathbb{Z}[P]/\langle p^{2}-p\colon p\in P\rangle$, \cite[\S
3]{rota-measure-logical}. This result was called
\emph{characterization theorem}. Here, we will generalize that
result for arbitrary bounded lattices.

In \cite{geis1}, L. Geissinger defines a ring, denoted
$\widetilde{V}(X)$, associated to a general lattice $X$. It is the
quotient between the free module $\mathbb{Z}^{\oplus X}$ with the
product $x\cdot y=x\wedge y$, and the ideal $\widetilde{J}$
generated by
\[
\{x\vee y+x\wedge y-x-y\,\colon\,x,y\in X\}.
\]
Geissinger characterizes the ring $\widetilde{V}(X)$ as representing
\emph{distributive valuations} (maps satisfying the
inclusion-exclusion principle for $k=2,3$). In \cite[p.466]{geis1},
Geissinger also noticed that any distributive valuation must satisfy
the inclusion exclusion principle for all $k\ge 2$.

\subsection*{Main results}

Let $X$ be a bounded lattice. The ring $\mathbb{Z}[X]/I$ is called
the \emph{lattice ring} of $X$, where $\mathbb{Z}[X]$ is the
polynomial ring generated by $X$ and $I=\langle 0_X,1_X-1,x\wedge
y-xy,x\vee y-x-y+xy\rangle$. In this article, all rings are assumed
to be commutative with an identity element. Morphisms between rings
respect addition, multiplication and the identity element.

Here, we provide a characterization of the generators of
$\widetilde{J}$ and the isomorphism between $\widetilde{V}(X)$ and
$\mathbb{Z}[X]/I$ given in Theorem \ref{repres}. We also extend the
results presented in \cite{geis1}, by proving the existence of a
Boolean lattice $Y$ such that $\mathbb{Z}[X]/I\cong\mathbb{Z}[Y]/I$
given in Theorem \ref{teo-Y}.

First, we prove the representability of the functor of bounded
measures:

\begin{theorem}
Let $X$ be a bounded lattice where a group $G$ acts.
Let $A$ be a ring
and let $B$ be an $A$-module
where $G$ acts trivially. Then,
\[
\mathcal{M}(X,B)\cong\Hom_\mathbb{Z}(\mathbb{Z}[X]/I,B)
\cong\Hom_A(A[X]/I_A,B),
\]
\[
\mathcal{M}(X,B;G)\cong\Hom_\mathbb{Z}(\mathbb{Z}[X]/I,B)^G
\cong\Hom_A(A[X]/I_A,B)^G,
\]
where $\mathcal{M}(X,B)$ is the space of $B$-valued bounded measures
on $X$, $\mathcal{M}(X,B;G)$ the space of $G$-invariant $B$-valued
measures on $X$ and $I_A$ is the ideal extended to $A[X]$.
\end{theorem}

Our second result is to prove that measures over bounded lattices
are in bijection with certain measures over Boolean lattices:

\begin{theorem}
Let $X$ be a bounded lattice. There exists a Boolean lattice $Y$
such that the functor $\mathcal{M}(X,-)$ is naturally isomorphic to
$\mathcal{M}(Y,-)$. In particular,
$\mathbb{Z}[X]/I\cong\mathbb{Z}[Y]/I$ and the study of measures over
a bounded lattice is equivalent to the study of measures over a
Boolean lattice. The lattice $Y$ is the Boolean lattice of
closed-open subsets on
$\mathcal{X}_{\mathbb{F}_2}=\Spec(\mathbb{F}_2[X]/I_{\mathbb{F}_2})$.
The ring $\mathbb{F}_2[X]/I_{\mathbb{F}_2}$ is reduced,
0-dimensional and its spectrum $\mathcal{X}_{\mathbb{F}_2}$ is a
Stone space.
\end{theorem}

As a third result we describe the ring $\mathbb{Z}[X]/I$ in the Noetherian case:

\begin{theorem}
Let $X$ be a bounded lattice. If
$\mathbb{Z}[X]/I$ is Noetherian (this is the case, for example, when
$X$ is finite), then there exists $0\leq n<\infty$ such that
\[
\mathbb{Z}[X]/I\cong\mathbb{Z}^n.
\]
\end{theorem}

Regarding finite lattices, we prove that \emph{measurability}
completely characterizes measures. The \emph{measurability} of a
bounded lattice $X$ is defined as $n(X)=\rk(\mathbb{Z}[X]/I)$:

\begin{theorem}
Let $X$ be a finite lattice. Let $A$ be a ring and let $B$ be an
$A$-module. Then, $\mathcal{M}(X,B)$ is isomorphic (as $A$-module)
to $B^n$, where $n=n(X)<\infty$. If the number $n$ is zero, the
lattice is called \emph{non-measurable} (as for example, the lattice
$M3$).
\end{theorem}

Summing up, our approach gives a unifying treatment to measures over
bounded lattices, by generalizing the Characterization Theorem of
Rota and the definition of valuation rings to any bounded lattice.
We think that working with bounded lattices is a relevant
generalization. As an example, we notice that, in many situations,
when a group acts on a lattice, it is in general false that the
orbit lattice $X/G$ is Boolean (not even distributive nor modular).
More examples are provided in Section \ref{s:Examples}.

Research on measure theory spreads on different areas, such as
functional analysis, mathematical physics, theory of ordered
structures and quantum probability theory. Here, we formulate the
problem of studying bounded measures in such a way to provide a
closer connection to algebraic geometry. This is done by defining a
functor of measures and asking for its representability (this is a
standard idea in deformation theory \cite{MR0146040} and category
theory \cite{MR1712872}). More specifically, we prove that given a
bounded lattice $X$, there exists a scheme $\mathcal{X}$ such that
measures on $X$ are in bijection with (scheme-theoretic) measures on
$\mathcal{X}$ (see Definition \ref{st-meas}). While the different
approaches used to study measures over lattices do not rely on
algebraic geometry, it is important to remark that the techniques
presented in \cite{MR2595208,MR2351269} are similar to the ones used
in this article. The main difference is that we work with more
general lattices than distributive ones. It is also important to
remark that our method can be easily extended to the study of
\emph{invariant} measures (that is, measures which are stable under
a given group of transformations).

\subsection*{Summary}

In Section \ref{s:Gen}, we give several elementary definitions. In
Section \ref{s:Meas}, we define the Measure functor $\mathcal{M}(X,-)$
and prove its representability by an Abelian group $M(X)$. We prove
that $M(X)$ has a natural structure of ring, and we characterize
this structure as a quotient of the polynomial ring generated by
$X$, $\mathbb{Z}[X]/I$. Next, we prove several universal properties
of this ring. We show that any measure on a bounded lattice can be
given as a measure on a Boolean lattice. We describe this Boolean
lattice as the lattice of closed-open subsets on a scheme
$\mathcal{X}_{\mathbb{F}_2}=\Spec(\mathbb{F}_2[X]/I_{\mathbb{F}_2})$.
In Section \ref{s:Examples}, we show relevant examples of measures over
lattices. In Section \ref{s:Finite}, we describe measures on finite
lattices by defining a new invariant, the \emph{measurability} (the
rank of $\mathbb{Z}[X]/I$). We compute this number for every lattice
with six elements or less.

\section{General definitions}\label{s:Gen}

The following general definitions can be
found in many books about lattice
theory, for example \cite{MR598630,MR931777}
\begin{itemize}
\item A poset $(X,\leq)$ is called a \emph{lattice} if every two elements $x,y$ have a supremum $x\vee y$
and an infimum $x\wedge y$. The operations $\vee$ and $\wedge$ satisfy associativity, commutativity,
idempotency and absorption
$$x\vee (y\vee z)=(x\vee y)\vee z,\quad x\vee y=y\vee x,\quad x\vee x=x,$$
$$x\wedge (y\wedge z)=(x\wedge y)\wedge z,\quad x\wedge y=y\wedge x,\quad x\wedge x=x,$$
$$x\wedge (x\vee y)=x,\quad x\vee (x\wedge y)=x.$$
We denote as $X$ to the lattice $(X,\leq,\vee,\wedge)$.

\item A lattice $X$ is called \emph{bounded}
if there exist two elements $0_X,1_X\in X$ such that
$$0\wedge x=0,\quad 0\vee x=x,\quad 1\wedge x=x,\quad 1\vee x=1.$$

\item A morphism (or a map or a function)
between bounded lattices $f:X\rightarrow Y$ is a
poset function such that
\[
f(x\wedge_X y)=f(x)\wedge_Y f(y),\quad f(x\vee_X y)=f(x)\vee_Y f(y),
\]
\[
f(0_X)=0_Y,\quad  f(1_X)=1_Y.
\]

\item A lattice $X$ is called \emph{complemented} if
it is bounded and for every $x\in X$ there exists $y\in X$ such that
$$x\vee y=1,\quad x\wedge y=0.$$

\item A lattice $X$ is called \emph{distributive} if for every $x,y,z\in X$,
\[
x\vee(y\wedge z)=(x\vee y)\wedge (x\vee z),\quad
x\wedge(y\vee z)=(x\wedge y)\vee (x\wedge z).
\]
In a distributive lattice, complements are unique, \cite[Ch.1,
Birkhoff Th.]{MR931777}. Hence, if $x$ has a complement, we denote
it as $x^\bot$.

\item A lattice is called \emph{Boolean} if it is
complemented and distributive.
Any distributive lattice can be extended
in a unique way to a Boolean lattice, \cite[Th. 13.18]{MR1501929}.

\item A lattice is called \emph{complete}
(resp. \emph{$\sigma$-lattice}) if any subset (resp. countable subset)
has both a join and a meet.
The \emph{join} and \emph{meet} of a subset $S$ are respectively the supremum (least upper bound) of $S$ 
and infimum (greatest lower bound) of $S$.

\item An element in a bounded lattice $z\in X$
is called an \emph{atom} if $z\neq 0$ and
\[
\{x\in X\,\colon\,0\leq x\leq z\}=\{0,z\}.
\]
\item A bounded lattice is called \emph{atomic} if
for every non-zero $x\in X$, there exists an atom $z$ such that $z\leq x$.
\end{itemize}

\section{Measures on bounded lattices as a functor}\label{s:Meas}

The measures defined below are valued in $A$-modules (where $A$ is a
ring). Thus, our definitions contain the usual notion of measure
(valued in the field of real numbers) as a particular case.

\begin{definition}\label{d:MeasureOnBoundedLattice}
Let $X$ be a bounded lattice, let $A$ be a ring and let $B$ be an
$A$-module. A function $\nu:X\rightarrow B$ is called a
\emph{measure} (or a \emph{bounded-measure}) with values in $B$ if
it satisfies
\[
\nu(0_X)=0_B,\quad
\nu(x_1\vee\dots\vee x_k)=
\sum_{s=1}^k(-1)^{s+1}\sum_{1\leq i_1<\dots<i_s\leq k}
\nu(x_{i_1}\wedge\dots\wedge x_{i_s})
\]
for every $x_1,\dots,x_k\in X$, $k\geq 2$. For example, if $k=2$,
\[
\nu(x\vee y)=\nu(x)+\nu(y)-\nu(x\wedge y),\quad\forall x,y\in X.
\]
If a group $G$ is acting on $X$
by lattice-automorphisms, an \emph{invariant measure} is a measure
stable under $G$.

The space of measures on $X$ with values in $B$
is an $A$-module denoted $\mathcal{M}(X,B)$
and the space of invariant measures
with values in $B$ is denoted as $\mathcal{M}(X,B;G)$.
By definition $\mathcal{M}(X,B)=\mathcal{M}(X,B;0)$,
where $0$ is the trivial group.
\end{definition}

\begin{proposition}
Let $X$ be a bounded lattice. Then,
there exists an Abelian group $M(X)$ representing
the functor $\mathcal{M}(X,-)$, i.e. we have
a natural isomorphism,
\[\Hom_\mathbb{Z}(M(X),-)\cong\mathcal{M}(X,-).\]
\end{proposition}
\begin{proof}
This proof follows easily from general principles.
Let $\mathbb{Z}^{\oplus X}$ be the free Abelian group
generated by $X$ and let $J$ be the submodule
generated by $0_X$ and
\[
x_1\vee\dots\vee x_k-
\sum_{s=1}^k(-1)^{s+1}\sum_{1\leq i_1<\dots<i_s\leq k}x_{i_1}\wedge\dots\wedge x_{i_s}
\]
for every $x_1,\dots,x_k\in X$, $k\geq 2$.
Let $M(X)=\mathbb{Z}^{\oplus X}/J$ and let $\pi:X\rightarrow M$ be the function given by $x\mapsto [x]$.

Let $B$ be an Abelian group and let
$\nu':M(X)\rightarrow B$ be a $\mathbb{Z}$-linear map.
Then, $\nu=\nu'\pi:X\rightarrow B$ is clearly a measure,
\begin{align*}
\nu(0_X)&=\nu'([0_X])=\nu'(0)=0_B,\\
\nu(x_1\vee\dots\vee x_k)&=\nu'([x_1\vee\dots\vee x_k])\\
&=
\nu'\left(\sum_{s=1}^k(-1)^{s+1}\sum_{1\leq i_1<\dots<i_s\leq k}
[x_{i_1}\wedge\dots\wedge x_{i_s}]
\right)\\
&=
\sum_{s=1}^k(-1)^{s+1}\sum_{1\leq i_1<\dots<i_s\leq k}
\nu'([x_{i_1}\wedge\dots\wedge x_{i_s}])\\
&=
\sum_{s=1}^k(-1)^{s+1}\sum_{1\leq i_1<\dots<i_s\leq k}
\nu(x_{i_1}\wedge\dots\wedge x_{i_s}).
\end{align*}
Also, $\nu'\pi\equiv 0$ implies $\nu'\equiv 0$. Then,
$\pi^*:\Hom_\mathbb{Z}(M(X),B)\rightarrow\mathcal{M}(X,B)$
is injective.

Now, a measure $\nu:X\rightarrow B$
extends to a linear
map $\nu':\mathbb{Z}^{\oplus X}\rightarrow B$
such that $\nu'(J)=0$,
\begin{align*}
\nu'\left(
x_1\vee\dots\vee x_k-
\sum_{s=1}^k(-1)^{s+1}\sum_{1\leq i_1<\dots<i_s\leq k}
x_{i_1}\wedge\dots\wedge x_{i_s}
\right)&=\\
\nu(x_1\vee\dots\vee x_k)-
\sum_{s=1}^k(-1)^{s+1}\sum_{1\leq i_1<\dots<i_s\leq k}
\nu(x_{i_1}\wedge\dots\wedge x_{i_s})
&=0.
\end{align*}
Then, $\nu'\in \Hom_\mathbb{Z}(M(X),B)$ and $\nu=\nu'\pi$.
Hence, $\pi^*$ is surjective.

Finally, if $f:B\rightarrow B'$ is a linear map,
the following diagram commutes,
\[
\xymatrix{
B\ar[d]_{f}&\Hom_\mathbb{Z}(M(X),B)\ar[r]^>>>>>{\pi^*}\ar[d]_{f_*}&\mathcal{M}(X,B)\ar[d]^{f_*}\\
B'&\Hom_\mathbb{Z}(M(X),B')\ar[r]^>>>>>{\pi^*}
\ar@{}[ur]|\circlearrowleft&\mathcal{M}(X,B')
}
\]
\end{proof}

\begin{remark}
In the next Theorem, we prove that
the wedge
product in $X$ induces
a ring structure on $M(X)$.
The fact that $X$ is bounded
implies that $M(X)$ becomes a
commutative ring with an identity element.

Let us define for a moment the notion of a non-normalized-measure.
Let $X$ be a
lattice, let $A$ be a ring and let $B$ be an $A$-module. We say the
$\mu:X\rightarrow B$ is a \emph{non-normalized-measure} if
\[
\mu(x_1\vee\dots\vee x_k)=
\sum_{s=1}^k(-1)^{s+1}\sum_{1\leq i_1<\dots<i_s\leq k}
\mu(x_{i_1}\wedge\dots\wedge x_{i_s})
\]
for every $x_1,\dots,x_k\in X$, $k\geq 2$.
Let $\mathcal{M}^{nn}(X,B)$ be the $A$-module of non-normalized measures
with values in $B$. Notice that $\mathcal{M}^{nn}(X,B)$
contains the constant functions.

A similar construction as before proves that the functor
$\mathcal{M}^{nn}(X,-)$ is representable by an Abelian group
$M^{nn}(X)$. Furthermore, the wedge product in $X$ induces a
structure of a commutative ring without an identity element in
$M^{nn}(X)$. In this article, we restrict to the case of
bounded-measures over a bounded lattice due to the fact that $M(X)$
becomes a commutative ring with an identity element. Although we
have chosen to study $\mathcal{M}(X,B)$, the analysis of
$\mathcal{M}^{nn}(X,B)$ is not relegated; the spaces
$\mathcal{M}^{nn}(X,B)$ and $\mathcal{M}(X,B)$ are related. For
example, if $X$ is bounded, then the map $\mu\mapsto
(\mu-\mu(0),\mu(0))$ is an $A$-module isomorphism between
$\mathcal{M}^{nn}(X,B)$ and $\mathcal{M}(X,B)\times B$.
\end{remark}

\begin{definition}
Let $X$ be a bounded lattice and let $\mathbb{Z}[X]$
be the ring generated by $X$,
that is, formal polynomials in the elements of $X$.
Consider the ideal of $\mathbb{Z}[X]$ generated by
\[
I=\langle
0_X,\quad 1_X-1,\quad x\wedge y-xy,\quad x\vee y-x-y+xy
\rangle.
\]
We call $\mathbb{Z}[X]/I$ the \emph{lattice ring} of $X$.

In the quotient ring $\mathbb{Z}[X]/I$ the element $0\in\mathbb{Z}$
and the class of $0\in X$ coincide (the same holds for $1$). Also
the meet operation is equal to the product operation and the join
operation satisfies
\[
x_1\vee\dots\vee x_n = 1-(1-x_1)\dots(1-x_n),\quad \forall x_1,\dots,x_n\in X.
\]
Indeed, by induction, if $n=2$, then
\[
1-x_1\vee x_2=1-x_1-x_2+x_1x_2=(1-x_1)(1-x_2),\quad \forall x_1,x_2\in X.
\]
If $n>2$, then
\begin{align*}
1-x_1\vee(x_2\vee\dots\vee x_n)&=(1-x_1)(1-x_2\vee\dots\vee x_n)\\
&=(1-x_1)(1-x_2)\dots(1-x_n),\quad \forall x_1,\dots,x_n\in X.
\end{align*}

In $\mathbb{Z}[X]/I$, complements are identified: if $y$ and $y'$
are two complements of $x$, then $x+y=1=x+y'$, hence, $y=y'$. Also,
the distributive law is formally satisfied,
\[
(x\vee y)\wedge (x\vee z)=(x+y-xy)(x+z-xz)=x+yz-xyz=x\vee(y\wedge z).
\]
\[
(x\wedge y)\vee (x\wedge z)=xy+xz-xyz=x(y+z-yz)=x\wedge(y\vee z).
\]
We denote $\mathcal{X}=\Spec(\mathbb{Z}[X]/I)$ to the spectrum
of $\mathbb{Z}[X]/I$. Regular functions on $\mathcal{X}$ are the
same as elements in $\mathbb{Z}[X]/I$ and closed points of
$\mathcal{X}$ correspond to maximal ideals.
\end{definition}

\begin{theorem}\label{repres}
Let $X$ be a bounded lattice. Then,
the Abelian group $M(X)$ representing the functor
$\mathcal{M}(X,-)$ has a structure of ring and as a ring
it is isomorphic to $\mathbb{Z}[X]/I$.
The ring structure on $M(X)$ is functorial in $X$.
\end{theorem}
\begin{proof}
Recall that $M(X)=\mathbb{Z}^{\oplus X}/J$, where
$J$ is generated by $0_X$ and
\[
x_1\vee\dots\vee x_k-
\sum_{s=1}^k(-1)^{s+1}\sum_{1\leq i_1<\dots<i_s\leq k}x_{i_1}\wedge\dots\wedge x_{i_s}
\]
for every $x_1,\dots,x_k\in X$, $k\geq 2$. In this proof we use
another notation for the generators:

\[
x_1\vee\dots\vee x_k+
\sum_{S\subseteq\{x_1,\dots,x_k\}}(-1)^{\# S}\bigwedge S.
\]

Let us first define a product in $\mathbb{Z}^{\oplus X}$
as $x\cdot y:=x\wedge y$.
By the lattice axioms, this product is associative, commutative
and has an identity element $1_X\in X$.
The zero element in $\mathbb{Z}^{\oplus X}$
is different from the vector $0_X$.

Let us prove that $J$ becomes an ideal with this product.
First, $x\cdot 0_X=x\wedge 0_X=0_X\in J$ for all $x\in X$.
Second, let $z\in X$ and consider the expressions $m$, $m_1$ and $m_2$,
\begin{align*}
m&:=z\cdot\left(x_1\vee\dots\vee x_k-
\sum_{s=1}^k(-1)^{s+1}\sum_{1\leq i_1<\dots<i_s\leq k}
x_{i_1}\wedge\dots\wedge x_{i_s}\right),\\
m_1&:=z\vee(x_1\vee\dots\vee x_k)-z-x_1\vee\dots\vee x_k+
z\wedge(x_1\vee\dots\vee x_k),\\
m_2&:=z\vee x_1\vee\dots\vee x_k+
\sum_{S\subseteq\{z,x_1,\dots,x_k\}}(-1)^{\#S}\bigwedge S\\
&=
z\vee x_1\vee\dots\vee x_k\\
&+
\sum_{S\subseteq\{x_1,\dots,x_k\}}(-1)^{\#S+1}z\wedge \bigwedge S+
\sum_{S\subseteq\{x_1,\dots,x_k\}}(-1)^{\#S}
\bigwedge S.
\end{align*}
Note that the expression in the middle,
when $S=\emptyset$, is equal to $-z$.
In other words,
\[
\sum_{S\subseteq\{x_1,\dots,x_k\}}(-1)^{\#S+1}z\wedge \bigwedge S=
-z+\sum_{s=1}^k(-1)^{s+1}\sum_{1\leq i_1<\dots<i_s\leq k}
z\wedge x_{i_1}\wedge\dots\wedge x_{i_s}.
\]
Then,
\begin{align*}
m-m_1+m_2&=
-\sum_{s=1}^k(-1)^{s+1}\sum_{1\leq i_1<\dots<i_s\leq k}z\wedge x_{i_1}\wedge\dots\wedge x_{i_s}\\
&+
z-z\vee(x_1\vee\dots\vee x_k)+x_1\vee\dots\vee x_k+m_2\\
&=
-\sum_{S\subseteq\{x_1,\dots,x_k\}}(-1)^{\#S+1}z\wedge \bigwedge S
-z\vee (x_1\vee\dots\vee x_k)\\
&+
x_1\vee\dots\vee x_k+m_2\\
&=
x_1\vee\dots\vee x_k
+\sum_{S\subseteq\{x_1,\dots,x_k\}}(-1)^{\#S}\bigwedge S.
\end{align*}
Given that $m_1,m_2$ and $m-m_1+m_2$ are in $J$, we deduce
that $m$ is also in $J$. Then, $M(X)$ has a ring structure.
Let us check the functoriality.
A lattice map $f:X\rightarrow Y$
induces a linear map $f_*:\mathbb{Z}^{\oplus X}\rightarrow
\mathbb{Z}^{\oplus Y}$
compatible with the ring structure
and $f_*(J_X)\subseteq J_Y$.
Then, $f_*$ induces a ring map $f_*:M(X)\rightarrow M(Y)$
such that the following diagram commutes,
\[
\xymatrix{
X\ar[d]_f & M(X)\otimes_{\mathbb{Z}}M(X)\ar[r]^-{\wedge}\ar[d]_{f_*\otimes f_*}\ar@{}[dr]|\circlearrowleft& M(X)\ar[d]^{f_*}\\
Y&M(Y)\otimes_{\mathbb{Z}}M(Y)\ar[r]^-{\wedge}& M(Y)
}
\]

It remains to check the isomorphism $M(X)\cong \mathbb{Z}[X]/I$.
Consider the surjective ring map
$\phi:\mathbb{Z}[X]\rightarrow M(X)$
and the surjective group map $\psi:\mathbb{Z}^{\oplus X}\rightarrow \mathbb{Z}[X]/I$. Both maps induced by the identity
$Id_X:X\rightarrow X$.
Note that $\psi$ is also a ring map
\[
\psi(xy)=\psi(x\wedge y)=[x\wedge y]=[x][y]=\psi(x)\psi(y)
\]
and $\phi(I)=0$ and
$\psi(J)=0$. Let us check $\psi(J)=0$,
\begin{align*}
\psi(x_1\vee\dots\vee x_k)&=[x_1\vee\dots\vee x_k]\\
&=1-(1-[x_1])\dots(1-[x_k])\\
&=
\sum_{s=1}^k(-1)^{s+1}\sum_{1\leq i_1<\dots<i_s\leq k}[x_{i_1}]\dots[x_{i_s}]\\
&=
\sum_{s=1}^k(-1)^{s+1}\sum_{1\leq i_1<\dots<i_s\leq k}[x_{i_1}\wedge\dots\wedge x_{i_s}].
\end{align*}
Then, $\psi(J)=0$. Finally, it is easy to
check $\psi\phi=Id_{\mathbb{Z}[X]/I}$ and $\phi\psi=Id_{M(X)}$.
\end{proof}

\begin{remark}[Relation with Geissinger's construction]
In \cite{geis1}, a ring is constructed, which is denoted by
$\widetilde{V}(X)$, and it is associated to a general lattice $X$.
It is the quotient between the free module $\mathbb{Z}^{\oplus X}$
with the product $x\cdot y=x\wedge y$ and the ideal $\widetilde{J}$
generated by
\[
\{x\vee y+x\wedge y-x-y\,\colon\,x,y\in X\}.
\]
From our previous developments, it is easy to prove the equality
$\widetilde{V}(X)/\langle 0_X\rangle=M(X)$. Indeed, denoting
\[
\iota_k(x_1,\dots,x_k):=x_1\vee\dots\vee x_k-
\sum_{s=1}^k(-1)^{s+1}\sum_{1\leq i_1<\dots<i_s\leq k}
x_{i_1}\wedge\dots\wedge x_{i_s},
\]
it follows from the proof of Theorem \ref{repres},
\[
z\cdot \iota_k(x_1,\dots,x_k)-\iota_2(z,x_1\vee\dots\vee x_k)+\iota_{k+1}(z,x_1,\dots,x_{k})
=
\iota_k(x_1,\dots,x_k).
\]
Given that $\widetilde{J}$ is an ideal,
if $\iota_k(x_1,\dots,x_k)\in \widetilde{J}$,
then $z\cdot  \iota_k(x_1,\dots,x_k)\in \widetilde{J}$. Then,
from the previous equation $\iota_{k+1}(z,x_1,\dots,x_{k})\in \widetilde{J}$ for all $z,x_1,\dots,x_k\in X$.
In particular, $J=\widetilde{J}+\langle 0_X\rangle$.
\end{remark}

\begin{notation}
Let $G$ be a group and let
$Z$ be an Abelian group
on which $G$ acts by additive maps on the left. We write $gz$ for the
action of $g\in G$ on $z\in Z$.
We denote $Z^G$ to the \emph{invariant subgroup} of $Z$
and $Z_G$ to the \emph{coinvariants} of $Z$,
\[
Z^G:=\{z\in Z\,\colon\,gz=z\quad\forall g\in G\},
\]
\[
Z_G:=Z/\text{submodule generated by }\{gz-z\,\colon\,g\in G,\,z\in Z\},
\]
\end{notation}

\begin{corollary}\label{univ-prop}
Let $X$ be a bounded lattice where a group $G$ acts.
Let $B$ be an Abelian group where $G$ acts trivially.

The map $\pi:X\rightarrow \mathbb{Z}[X]/I$ is a measure and satisfies the following universal property.
Any measure with values in $B$ factorizes as
a $\mathbb{Z}$-linear map $\mathbb{Z}[X]/I\rightarrow B$.
Also, the map $\pi_G:X\rightarrow (\mathbb{Z}[X]/I)_G$ is an invariant measure
and satisfies the following universal property.
Any invariant measure with
values in $B$ factorizes as a $\mathbb{Z}$-linear
map $(\mathbb{Z}[X]/I)_G\rightarrow B$.

We can represent the properties of $\pi$ and $\pi_G$
by the following diagrams,
\[
\xymatrix{
X\ar[r]^{\forall\nu}\ar[d]_{\pi}
\ar@{}[dr]|<<<<\circlearrowleft&B\\
\mathbb{Z}[X]/I\ar@{-->}[ur]_{\exists!\overline{\nu}}&
}\qquad
\xymatrix{
X\ar[r]^{\forall\nu'}\ar[d]_{\pi_G}
\ar@{}[dr]|<<<<\circlearrowleft&B\\
(\mathbb{Z}[X]/I)_G\ar@{-->}[ur]_{\exists!\overline{\nu'}}&
}
\]
where $\nu$ (resp. $\nu'$) is a measure (resp. invariant measure)
and $\overline{\nu}$, $\overline{\nu'}$ are linear maps. The
commutativity means that $\nu=\overline{\nu}\pi$,
$\nu'=\overline{\nu'}\pi_G$. Alternatively,
\[
\pi^*:\Hom_\mathbb{Z}(\mathbb{Z}[X]/I,B)\cong\mathcal{M}(X,B)
\quad\text{and}\quad
\]
\[
\pi_G^*:\Hom_\mathbb{Z}((\mathbb{Z}[X]/I)_G,B)\cong\mathcal{M}(X,B;G).
\]
\end{corollary}
\begin{proof}
From Theorem \ref{repres}, it follows that
$\pi$ is a measure with values in $\mathbb{Z}[X]/I$ and satisfies the universal property.
Let us define $\pi_G$ and prove its universal property.

Let $\nu':X\to B$ be an invariant measure. From the universal property of $\pi$
there exists a linear map $\nu'':\mathbb{Z}[X]/I\to B$
such that $\nu'=\nu''\pi$. Clearly, $\nu''$
is $G$-invariant.
Let $p$ be the quotient $p:\mathbb{Z}[X]/I\to(\mathbb{Z}[X]/I)_G$
and consider the following diagram,
\[
\xymatrix{
X\ar[r]^{\nu'}\ar[d]_{\pi}
\ar@{}[dr]|<<<<<\circlearrowleft
\ar@{}[dr]|>>>>>\circlearrowleft
&B\\
\mathbb{Z}[X]/I\ar@{-->}[ur]^{\nu''}\ar[r]^{p}&(\mathbb{Z}[X]/I)_G
\ar@{-->}[u]_{\exists!\overline{\nu'}}
}
\]
Recall that any $G$-linear map into a trivial $G$-module factorizes uniquely
over a map from $(\mathbb{Z}[X]/I)_G$.
Then, there exists
a unique map $\overline{\nu'}:(\mathbb{Z}[X]/I)_G\to B$
such that $\nu''=\overline{\nu'}p$.
The result follows by defining $\pi_G:=p\pi$ and noting
$\nu'=\nu''\pi=\overline{\nu'}p\pi=\overline{\nu'}\pi_G$.
\end{proof}

\begin{definition}
We call $\pi$ the \emph{universal measure} and $\pi_G$ the
\emph{universal invariant measure}. If no confusion arises, we write
$x$ instead of $\pi(x)$ (or $\pi_G(x)$).
\end{definition}

\begin{corollary}
Let $X$ be a bounded lattice where a group $G$ acts.
Let $A$ be a ring
and let $B$ be an $A$-module
where $G$ acts trivially.

Then, the universal measure and the universal invariant
measure induce the following
three equivalent characterization of
the spaces $\mathcal{M}(X,B)$ and $\mathcal{M}(X,B;G)$,
\[
\mathcal{M}(X,B)\cong\Hom_\mathbb{Z}(\mathbb{Z}[X]/I,B)
\cong\Hom_A(A[X]/I_A,B),
\]
\[
\mathcal{M}(X,B;G)\cong\Hom_\mathbb{Z}(\mathbb{Z}[X]/I,B)^G
\cong\Hom_A(A[X]/I_A,B)^G.
\]
\end{corollary}
\begin{proof}
First, from the $\otimes$-$\Hom$ adjunction we have,
\[
\Hom_A(A[X]/I_A,B)\cong
\Hom_A(\mathbb{Z}[X]/I\otimes_\mathbb{Z}A,B)
\cong
\]
\[
\Hom_\mathbb{Z}(\mathbb{Z}[X]/I,\Hom_A(A,B))\cong
\Hom_\mathbb{Z}(\mathbb{Z}[X]/I,B).
\]
Also, recall that the functor $(-)_G$ is naturally isomorphic to $(-)\otimes_{\mathbb{Z}[G]}\mathbb{Z}$
and the functor $(-)^G$ is naturally isomorphic to $\Hom_{\mathbb{Z}[G]}(\mathbb{Z},-)$, \cite[Lemma 6.1.1]{MR1269324}.
Then, using again the $\otimes$-$\Hom$ adjunction, we have
\[
\Hom_{\mathbb{Z}}((\mathbb{Z}[X]/I)_G,B)\cong
\Hom_{\mathbb{Z}}(\mathbb{Z}\otimes_{\mathbb{Z}[G]}\mathbb{Z}[X]/I,B)\cong
\]
\[
\Hom_{\mathbb{Z}[G]}(\mathbb{Z},\Hom_\mathbb{Z}(\mathbb{Z}[X]/I,B))\cong
\Hom_\mathbb{Z}(\mathbb{Z}[X]/I,B)^G.
\]
\end{proof}

\begin{definition}\label{st-meas}
Let $A$ be a ring and let $B$ be a $A$-module.
Consider the affine scheme
$\mathcal{X}=\Spec(A[X]/I_A)$ over
$\mathcal{S}=\Spec(A)$ with canonical map $f:\mathcal{X}\to\mathcal{S}$
and let $\widehat{B}$ be the sheaf of $\mathcal{O}_{\mathcal{S}}$-modules
associated to $B$.
Following \cite[III.7]{MR2018901}, we can adapt the definition
of a measure to the context of schemes.
We define the \emph{sheaf of measures over the lattice $X$}
as the sheaf in $\mathcal{S}$ of $\mathcal{O}_{\mathcal{S}}$-modules
\[
\mathcal{H}om_{\mathcal{O}_{\mathcal{S}}}(f_{*}\mathcal{O}_{\mathcal{X}},\mathcal{O}_{\mathcal{S}})
\]
and we define the sheaf of
\emph{measures over the lattice $X$ with values in $B$} as
\[
\mathcal{F}_{B}:=\mathcal{H}om_{\mathcal{O}_{\mathcal{S}}}(f_*\mathcal{O}_{\mathcal{X}},\widehat{B}).
\]
Since $\Gamma(\mathcal{S},f_*\mathcal{O}_{\mathcal{X}})=
\Gamma(\mathcal{X},\mathcal{O}_{\mathcal{X}})=A[X]/I_A$, the space of global sections of $\mathcal{F}_{B}$ is
\[
\Gamma(\mathcal{S},\mathcal{F}_{B})\cong
\Hom_A(\Gamma(\mathcal{X},\mathcal{O}_{\mathcal{X}}),B)=
\Hom_A(A[X]/I_A,B).
\]
Then, from the previous Corollary, we get
\[
\mathcal{M}(X,B)\cong
\Hom_A(\Gamma(\mathcal{X},\mathcal{O}_{\mathcal{X}}),B)
,\qquad
\]
\[
\mathcal{M}(X,B;G)\cong
\Hom_A(\Gamma(\mathcal{X},\mathcal{O}_{\mathcal{X}}),B)^G.
\]
\end{definition}

The next result generalizes
the functor given in
\cite{rota-valuation-ring}
to bounded lattices

\begin{corollary}
Let $A$ be a ring,
let \textbf{A-Mod} be the category of $A$-modules
and let \textbf{Bou} be the category of
bounded lattices.
Then,
we have a functor
\[
\mathcal{M}(-,-):\textbf{Bou}^{op}\times
\textbf{A-Mod}\rightarrow \textbf{A-Mod}.
\]
\end{corollary}
\begin{proof}
Let $X$ be a bounded lattice and let $B$ be an $A$-module.
Clearly $\mathcal{M}(X,B)=\Hom_A(A[X]/I_A,B)$ is an $A$-module.
Let $f:X\rightarrow X'$ be a map between two bounded
lattices and let $g:B\rightarrow B'$ be a map between two
$A$-modules. Then $(f^*,g_*):\mathcal{M}(X',B)\rightarrow
\mathcal{M}(X,B')$, $\nu\mapsto g\nu f$, is $A$-linear.
\end{proof}

\begin{theorem}\label{teo-Y}
Let $X$ be a bounded lattice
and consider the natural map $\pi:X\to\mathbb{Z}[X]/I$.
Let $Y$ be the Boolean lattice generated by the image of $\pi$.
Then,
\[
\pi^*:\mathcal{M}(Y,-)\rightarrow\mathcal{M}(X,-)
\]
is
a natural isomorphism.
In particular, $\mathbb{Z}[X]/I\cong\mathbb{Z}[Y]/I$ and
the study of measures over bounded lattices is equivalent to the
study of measures over Boolean lattices.
\end{theorem}
\begin{proof}
Recall that for any commutative ring with unity, the set $E$ of idempotents
forms a Boolean lattice with $x\wedge y := xy$ and $x\vee y := x+y-xy$. Let us prove
that the Boolean lattice $Y\subseteq E$ generated by the image of $\pi:X\to\mathbb{Z}[X]/I$ satisfies,
\[
\pi^*:\mathcal{M}(Y,-)\cong \mathcal{M}(X,-).
\]
Indeed, if $B$ is an Abelian group and $\nu:Y\to B$ a measure such that $\nu\pi:X\to B$ is zero,
then, $\nu$ is zero over the image of $\pi$ (which is a distributive lattice) and by
\cite[Cor. 2.2.2]{MR1670882}, $\nu$ extends uniquely to the Boolean lattice $Y$.
The uniqueness of the extension implies that it must be zero. Hence, $\pi^*$ is injective.

Let us prove that
$\pi^*:\mathcal{M}(Y,B)\rightarrow\mathcal{M}(X,B)$ is surjective.
Let $\nu'$ be a measure in $X$.
Define $\nu$ over the image of $\pi$ as $\nu(\pi(x)):=\nu'(x)$ for $x\in X$.
Let us see that $\nu$ is well-defined.
Given that $\nu'$ factorizes through $\pi$, if
$\pi(x_1)=\pi(x_2)$, then $\nu'(x_1)=\nu'(x_2)$.
This implies that $\nu$ does not depend on the representative of $\pi(x)$. By \cite[Cor. 2.2.2]{MR1670882}, $\nu$
extends uniquely to $Y$. Hence, $\pi^*$ is surjective.

The naturality follows from the following diagram,
\[
\xymatrix{
X\ar[d]_f\ar[r]^\pi\ar@{}[dr]|\circlearrowleft&Y\ar[d]^{F}\\
X'\ar[r]^\pi&Y'
}
\]
where $F:Y\rightarrow Y'$ is defined as the restriction
of the natural ring map $f_{*}:\mathbb{Z}[X]/I\rightarrow \mathbb{Z}[X']/I$.
Indeed, from the commutativity of the previous diagram, we obtain
the following commutativity
\[
\xymatrix{
X\ar[d]_f&
\mathcal{M}(Y',-)\ar[d]_{F^*}\ar[r]^{\pi^*}\ar@{}[dr]|\circlearrowleft&\mathcal{M}(X',-)\ar[d]^{f^{*}}\\
X'&\mathcal{M}(Y,-)\ar[r]^{\pi^*}&\mathcal{M}(X,-)
}
\]
\[
\pi^{*} F^{*}(\nu)=\nu F\pi=\nu \pi f= f^*\pi^*(\nu).
\]

Finally, from the naturality we deduce the isomorphism of groups $\pi_*:M(X)\rightarrow M(Y)$
and it is easy to check that $\pi_*$ is compatible with the ring structures. Then,
we get a natural isomorphism of rings
$\mathbb{Z}[X]/I\cong\mathbb{Z}[Y]/I$.
\end{proof}

\begin{corollary}\label{car-Y}
Let $X$ be a bounded lattice.
Then, the Boolean lattice $\Sigma_{\mathbb{F}_2}$,
of closed and open subsets on
$\mathcal{X}_{\mathbb{F}_2}=\Spec(\mathbb{F}_2[X]/I_{\mathbb{F}_2})$
is isomorphic to the Boolean lattice $Y$ of Theorem \ref{teo-Y}.
\end{corollary}
\begin{proof}
The isomorphism $Y\cong \mathbb{F}_2[Y]/I_{\mathbb{F}_2}$
can be found in \cite[Cor. p. 234]{rota-measure-logical} or below in Proposition \ref{bool-f2}.
The isomorphism
$\mathbb{F}_2[X]/I_{\mathbb{F}_2}\cong \Sigma_{\mathbb{F}_2}$ follows
from Stone's Theorem, \cite[13.7]{kleiman2013term}
and
$\mathbb{F}_2[Y]/I_{\mathbb{F}_2}\cong\mathbb{F}_2[X]/I_{\mathbb{F}_2}$ follows
from Theorem \ref{teo-Y}.
\end{proof}

\section{Some examples}\label{s:Examples}

\begin{example}
An interesting example is
the space of measures with values in $\mathbb{F}_2$.
In this case, the ring $\mathbb{F}_2[X]/I_{\mathbb{F}_2}$ is Boolean
(every element is idempotent)
and by Stone's Theorem, \cite[13.7]{kleiman2013term},
it is isomorphic to the space of continuous functions from
$\mathcal{X}_{\mathbb{F}_2}=\Spec(\mathbb{F}_2[X]/I_{\mathbb{F}_2})$ to $\mathbb{F}_2$,
\[
\mathbb{F}_2[X]/I_{\mathbb{F}_2}\cong\Cont(\mathcal{X}_{\mathbb{F}_2},\mathbb{F}_2).
\]
Hence,
$\mathcal{M}(X,\mathbb{F}_2)$ is
characterized as linear function from $\Cont(\mathcal{X}_{\mathbb{F}_2},\mathbb{F}_2)$
to $\mathbb{F}_2$.
\end{example}

It is possible to give to a Boolean lattice $Z$ a ring structure
with addition $z_1+z_2=(z_1\wedge z_2^\bot)\vee(z_1^\bot\wedge z_2)$
and multiplication $z_1z_2=z_1\wedge z_2$, \cite[Exerc. 24,
p.14]{MR0242802}. Also, we can associate to the lattice $Z$ the
Boolean ring $\mathbb{F}_2[Z]/I_{\mathbb{F}_2}$. The next result can
be found in \cite{rota-measure-logical}.

\begin{proposition}\label{bool-f2}
Let $Z$ be a Boolean lattice. Then, $Z$, viewed
as a ring, is isomorphic to $\mathbb{F}_2[Z]/I_{\mathbb{F}_2}$,
\[
Z\cong \mathbb{F}_2[Z]/I_{\mathbb{F}_2}.
\]
\end{proposition}
\begin{proof}
This proof is straightforward.
Let us prove that the universal measure with values in $\mathbb{F}_2$,
$\pi:Z\rightarrow\mathbb{F}_2[Z]/I_{\mathbb{F}_2}$ is an isomorphism of rings.
Notice the equality $\pi(z^\bot)=1-\pi(z)$. Then,
\begin{align*}
\pi(z_1z_2)&=\pi(z_1\wedge z_2)=\pi(z_1)\pi(z_2).\\
\pi(z_1+z_2)&=\pi((z_1\wedge z_2^\bot)\vee(z_1^\bot\wedge z_2))\\
&=
\pi(z_1\wedge z_2^\bot)+\pi(z_1^\bot\wedge z_2)-\pi(z_1\wedge z_2^\bot\wedge z_1^\bot\wedge z_2)\\
&=
\pi(z_1\wedge z_2^\bot)+\pi(z_1^\bot\wedge z_2)\\
&=
\pi(z_1)(1-\pi(z_2))+(1-\pi(z_1))\pi(z_2)\\
&=
\pi(z_1)+\pi(z_2)-2\pi(z_1)\pi(z_2)\\
&=
\pi(z_1)+\pi(z_2).
\end{align*}
The identity $Id:Z\rightarrow Z$ induces a ring map
$\mathbb{F}_2[Z]\rightarrow Z$
and by De Morgan's law, its kernel contains $I$.
For example,
\begin{align*}
x\vee y+x\wedge y &=
((x\vee y)\wedge (x\wedge y)^\bot)\vee((x\vee y)^\bot\wedge (x\wedge y))\\
&=
((x\vee y)\wedge (x^\bot\vee y^\bot))\vee
(x^\bot\wedge y^\bot\wedge x\wedge y)\\
&=
(x\vee y)\wedge (x^\bot\vee y^\bot)\\
&=
((x\vee y)\wedge x^\bot)\vee ((x\vee y)\wedge y^\bot)\\
&=
(y\wedge x^\bot)\vee (x\wedge y^\bot)\\
&=x+y.
\end{align*}
Hence, $\overline{Id}:\mathbb{F}_2[Z]/I_{\mathbb{F}_2}\rightarrow Z$ is well-defined
and inverse to $\pi$.
\end{proof}

\begin{example}
From the previous Proposition, we deduce that
there exists a lattice $X$
such that $\mathbb{Z}[X]/I$ is not Noetherian (resp. of finite type).
Let $X=\mathbb{F}_2[x_1,\dots]/\langle x_1^2-x_1,\dots\rangle$
be a Boolean ring with infinitely many variables viewed as a Boolean lattice.

Assume $\mathbb{Z}[X]/I$ is Noetherian (resp. of finite type), then
$\mathbb{F}_2[X]/I_{\mathbb{F}_2}\cong X$ is Noetherian (resp. of finite type),
a contradiction.
\end{example}

\begin{example}\label{powerset}
Consider the lattice $X=\mathcal{P}(\{x_1,\dots,x_n\})$ of subsets of a set with $n$ elements. The atoms
of this lattice are $\{x_1,\dots,x_n\}$ and every element in $\mathbb{Z}[X]/I$ is a polynomial in these variables.
Also, $x_ix_j=0$ and $x_i^2=x_i$, $i\neq j$. Then, the ring $\mathbb{Z}[X]/I$ is isomorphic to $\mathbb{Z}^n$
with addition and multiplications coordinate-wise.
Recall that a finite Boolean lattice is a power set, \cite[Ch.II.1, Cor.12]{MR509213}

In general, not every Boolean lattice is a power set.
It is known that
a Boolean lattice is isomorphic
to a power set if and only if it is complete and atomic,
see \cite[Ch.V,\S6,Th.18]{MR598630}.
The Boolean lattice of measurable subsets
in $\mathbb{R}^n$ is atomic but not complete
and the Boolean lattice of measurable subsets
in $\mathbb{R}^n$ modulo sets of measure zero
is complete but not atomic (it follows from
$|S|=|\overline{S}|=|S^\circ|$).
\end{example}

\begin{example}\label{ex:SubspacesHilbert}
Consider the lattice $X$ of subspaces in $\mathbb{R}^n$. The atoms
of this lattice are the one-dimensional subspaces. Every measure
invariant under the orthogonal group $O(n)$ is determined by its
value $\lambda$ at some one-dimensional subspace. Specifically, if
$\{v_1,\dots,v_k\}$ is a basis for $L$, then
\[
\nu(L)=\nu(v_1)+\dots+\nu(v_k)=\lambda\dim(L).
\]
Then, any $O(n)$-invariant measure on $X$ is, up to a constant, the
dimension.

A similar situation occurs with $Sym(S)$-invariant
measures on $\mathcal{P}(S)$, where $\# S<\infty$.
Any measure invariant under the group of symmetries on $S$
is, up to a constant, the cardinal.

We will discuss dimension functions in more general spaces in the
next example.
\end{example}

\begin{remark}\label{rem:DimensionFunctions}
In the theory of von Neumann algebras, dimension functions play a
key role in the classification of factors. Factors are classified by
appealing to the ranges of their dimension functions. As is well
known, the collection of projection operators of a factor can be
endowed with an orthomodular lattice structure.

In order to provide a definition of dimension function for abstract
orthomodular lattices, we need the notion of \textit{dimension
lattice}. Let ``$\equiv$" be an equivalence relation in a complete
orthomodular lattice $X$. A \textit{dimension lattice} is a pair
$(X,\equiv)$, where the following conditions hold
\cite{Kalmbach-OrthomodularLattices}:

\begin{itemize}
\item if $a\equiv 0$, then $a=0$.
\item If $a\perp b$ and $c\equiv a\vee b$, then there exist $d$ and
$e$ such that $a\equiv d$, $b\equiv e$, $d\perp e$ and $c=d\vee e$.
\item Let $A$ and $B$ be two orthogonal sets (i.e., all their elements are pairwise
orthogonal) with a bijection
$g:A\longrightarrow B$ such that $a\equiv g(a)$ for all $a$. Then,
we have $\bigvee A\equiv\bigvee B$.
\item Let $a\sim b$ if and only if there exists $c\in X$ such that $a\vee c=1_X=b\vee c$ and $a\wedge c=0_{X}=b\wedge c$. Then,
$a\sim b$ implies $a\equiv b$.
\end{itemize}

An element $x\in X$ is said to be finite whenever $x\equiv y$ and
$y\leq x$ imply $x=y$.
A real-valued dimension function $\alpha$ on a dimension lattice $X$
is a map from $X$ to the set $\mathbb{R}_{\geq 0}$ satisfying the
following conditions:

\begin{itemize}
\item $\alpha(0_{X})=0$
\item $\alpha(\bigvee A)=\sum_{a\in A}\alpha(a)$, for all non-empty
orthogonal and countable family in $X$.
\item $\alpha(a)<+\infty$ for all finite $a$.
\item $a\equiv b$ if and only if $\alpha(a)=\alpha(b)$.
\end{itemize}

As is well known, the set of projection operators
$\mathcal{P}(\mathcal{V})$ of a von Neumann algebra $\mathcal{V}$
can be endowed with a natural complete orthomodular lattice
structure. A key result is that $\mathcal{P}(\mathcal{V})$ generates
$\mathcal{V}$ in the sense that
$(\mathcal{P}(\mathcal{V}))''=\mathcal{V}$ (in words: the double
commutant of $\mathcal{P}(\mathcal{V})$ equals $\mathcal{V}$). Two
projections $a$ and $b$ in $\mathcal{V}$ are called equivalent (and
we denote $a\simeq b$) if there exists an operator in $\mathcal{V}$
that maps the vectors in $a^{\perp}$ into zero, and it is an
isometry between the subspaces of $a$ and $b$. In formulae: $a\simeq
b$ if there exists a partial isometry $v\in\mathcal{V}$ such that
$a=v^{\ast}v$ and $b=vv^{\ast}$. Clearly,
$\simeq$ is an equivalence relation
in $\mathcal{V}$ and we denote by
$\mathcal{P}(\mathcal{V})/{\simeq}$ to the set of its equivalence
classes. A partial order relation in $\mathcal{P}(\mathcal{V})$ can
be defined by appealing to $\simeq$: $a\preceq b$ if and only if
there exists $c$ such that $c\leq b$ and $a\simeq c$ (here ``$\leq$"
denotes the usual partial order between operators). Notice that
``$\simeq$" satisfies the properties of ``$\equiv$" for the
orthomodular lattice of projection operators in a factor von Neumann
algebra.

It can be proved that if $\mathcal{V}$ is a factor von Neumann
algebra, then $\preceq$ becomes a total order in
$\mathcal{P}(\mathcal{V})/{\simeq}$. This can be used to prove the
following: if $\mathcal{V}$ is a factor von Neumann algebra there
exists a dimension function
$d:\mathcal{P}(\mathcal{V})\longrightarrow[0,+\infty]$ with the
following properties:

\begin{enumerate}
\item $d(a)=0$ if and only if $a=0$.
\item If $a\perp b$, then $d(a\vee b)=d(a)+d(b)$.
\item $d(a)\leq d(b)$ if and only if $a \preceq b$.
\item $d(a)<+\infty$ if and only if $a$ is a finite projection (a projection $a$
is finite if for all $b$ such that $a\simeq b$ and $b\leq a$, it
follows that $a=b$).
\item $d(a)=d(b)$ if and only if $a\sim b$.
\item $d(a)+d(b)=d(a\vee b)+d(a\wedge b)$.
\end{enumerate}

Properties $1$ and $6$ guarantee that dimension functions in factor
von Neumann algebras are of the kind of measures on bounded lattices
studied in this paper (compare with definition
\ref{d:MeasureOnBoundedLattice}). Thus, we can apply all the
machinery of algebraic geometry introduced here for their study. In
particular, dimension functions can be seen as invariant measures
under the action of partial isometries. The above Example
\ref{ex:SubspacesHilbert} can be seen as a particular case of the
more general framework of dimension functions in orthomodular
lattices.
\end{remark}

\begin{example}
Let M3 be the following lattice,
\[
\xymatrix{
&\circ\ar@{-}[dl]\ar@{-}[d]\ar@{-}[dr]&\\
x_1\ar@{-}[dr]&x_2\ar@{-}[d]&x_3\ar@{-}[dl]\\
&\circ&
}
\]
The ring $\mathbb{Z}[\mathrm{M3}]/I$
has three generators, $x_1,x_2,x_3$
that satisfy, $x_1+x_2=x_1+x_3=x_2+x_3=x_1+x_2+x_3=1$.
Simplifying, we get $x_1=x_2=x_3=0$
and then $1=0$. Hence, $\mathbb{Z}[\mathrm{M3}]/I$ is the zero ring.

A more extreme case is the lattice $\mathrm{M}\omega$,
\[
\xymatrixrowsep{.3cm}
\xymatrixcolsep{.3cm}
\xymatrix{
&&\circ\ar@{-}[dll]\ar@{-}[dl]\ar@{-}[d]\ar@{-}[dr]\ar@{-}[drr]\\
\circ&\circ&\dots&\circ&\circ\\
&&\circ\ar@{-}[ull]\ar@{-}[ul]\ar@{-}[u]\ar@{-}[ur]\ar@{-}[urr]\\
}
\]
Notice that $\#(\mathrm{M}\omega)=\infty$ and also,
$\mathbb{Z}[\mathrm{M}\omega]/I=0$. The last example of this type is the following lattice $X$,
\[
\xymatrixrowsep{.3cm}
\xymatrixcolsep{.3cm}
\xymatrix{
&&\circ\ar@{-}[dl]\ar@{-}[ddrr]\\
&x_4\ar@{-}[dl]\ar@{-}[d]\ar@{-}[dr]\\
x_1\ar@{-}[dr]&x_2\ar@{-}[d]&x_3\ar@{-}[dl]&& y\\
&x_0\\
&&\circ\ar@{-}[ul]\ar@{-}[uurr]
}
\]
It is easy to prove $\mathbb{Z}[X]/I\cong\mathbb{Z}^2$
and the equalities $x_i=1-y$ in $\mathbb{Z}[X]/I$, $0\leq i\leq 4$.
\end{example}

\begin{example}
Let N5 and M2 be the following lattices,
\[
\xymatrixrowsep{.3cm}
\xymatrix{
&\circ\ar@{-}[dl]\ar@{-}[ddr]&\\
\circ\ar@{-}[dd]&&\\
&&\circ\ar@{-}[ddl]\\
\circ\ar@{-}[dr]&&\\
&\circ&
}\quad
\xymatrixrowsep{1cm}
\xymatrix{
&\circ\ar@{-}[dl]\ar@{-}[dr]&\\
\circ\ar@{-}[dr]&&\circ\ar@{-}[dl]\\
&\circ&
}
\]
The ring associated to N5 and to M2 are the same and equal to $\mathbb{Z}^2$. Then,
measures on N5 are the same as measures on M2,
\[
\mathcal{M}(\mathrm{N5},-)\cong\mathcal{M}(\mathrm{M2},-).
\]
A more extreme case is the following non-complete lattice N$\omega$,
\[
\xymatrixrowsep{.3cm}
\xymatrixcolsep{.6cm}
\xymatrix{
&\circ\ar@{-}[dl]\ar@{-}[ddr]\\
\ar@{.}[d]\\
\circ\ar@{-}[d]&&\circ\\
\circ\\
&\circ\ar@{-}[ul]\ar@{-}[uur]
}
\]
It is easy to prove that
the ring associated to N$\omega$ is also $\mathbb{Z}^2$.
\end{example}

\begin{proposition}\label{max-lat}
Let $X$ be the lattice
$0\leq 1\leq\dots\leq n$, where $i\vee j=\max(i,j)$ and $i\wedge j=\min(i,j)$. Then $\mathbb{Z}[X]/I\cong\mathbb{Z}^n$
as rings.
\end{proposition}
\begin{proof}
Let $\phi:\mathbb{Z}^{\oplus X}\rightarrow \mathbb{Z}^n$
be the surjective linear map
given by $\phi(0_X)=0$
and for $1\leq i\leq n$, $\phi(i)=e_1+\dots+e_i$, where $e_i$ is the vector
with $(e_i)_i=1$ and $(e_i)_j=0$ for $j\neq i$.
Let us prove that the ideal $J$ is equal to $J=\langle 0_X\rangle$,
\begin{align*}
\sum_{S\subseteq\{x_1\}}(-1)^{\# S}\bigwedge S&=0_X-x_1,\quad x_1\in X.\\
\sum_{S\subseteq\{x_1,x_2\}}(-1)^{\# S}\bigwedge S&=0_X-x_1-x_2+x_1\\
&=
0_X-x_2,\quad x_1\leq x_2\in X.
\end{align*}
For $k>2$, let $x_1\leq\dots\leq x_k$ in $X$,
\begin{align*}
\sum_{S\subseteq\{x_1,\dots,x_k\}}(-1)^{\# S}\bigwedge S
&=
\sum_{S\subseteq\{x_2,\dots,x_k\}}(-1)^{\# S+1}x_1\wedge\bigwedge S\\
&+
\sum_{S\subseteq\{x_2,\dots,x_k\}}(-1)^{\# S}\bigwedge S\\
&=
-x_1\left(\sum_{S\subseteq\{x_2,\dots,x_k\}}(-1)^{\# S}\right)+0_X-x_k\\
&=0_X-x_k.
\end{align*}
In the previous line we used the inductive hypothesis and the
following equation
\[
\sum_{S\subseteq\{x_2,\dots,x_k\}}(-1)^{\# S}=
\sum_{i=0}^{k-1}\sum_{\#S=i}(-1)^i=
\sum_{i=0}^{k-1}\binom{k-1}{i}(-1)^i=(1-1)^{k-1}=0.
\]
Then $J=\langle 0_X\rangle$. In particular,
$\phi$ induces an isomorphic linear map $M(X)\rightarrow \mathbb{Z}^n$.
Finally, assume $i<j$ and let us view $\mathbb{Z}^n$ as a ring.
Then,
\begin{align*}
\phi(i\wedge j)&=\phi(\min(i,j))\\
&=
\phi(i)\\
&=
e_1+\dots+e_i\\
&=
(e_1+\dots+e_i)(e_1+\dots+e_j)\\
&=
\phi(i)\phi(j).
\end{align*}
\end{proof}

\section{Measurability and measures on a finite lattice}\label{s:Finite}

\begin{remark}\label{bound}
Recall that
$\mathbb{Z}[X]/I$ is isomorphic, as a ring, to
$M(X)=\mathbb{Z}^{\oplus X}/J$.
In particular, if $X$ is finite
and $A$ is a field,
$A[X]/I_A$ is always Artinian (i.e. $\mathcal{X}_{A}$
is a finite union of points).
Furthermore, $\#\mathcal{X}_{A}\leq\#X-1$.
This bound is sharp, the lattice
$0\leq 1$ has $2$ elements and $\#\mathcal{X}_{A}=1$.
\end{remark}

\begin{proposition}\label{z-n}
Let $X$ be a bounded lattice.
If $\mathbb{Z}[X]/I$ is Noetherian or $X$ is finite,
then $\mathbb{Z}[X]/I$ is isomorphic
to $\mathbb{Z}^n$, where $0\leq n<\infty$. The case of the zero
ring is included with $n=0$.

Furthermore, if the rank of $\mathbb{Z}[X]/I$ is equal to $n<\infty$,
then $\mathbb{Z}[X]/I\cong\mathbb{Z}^n$.
\end{proposition}
\begin{proof}
By Theorem \ref{teo-Y} we can assume, without loss of generality,
that $X$ is a Boolean lattice.
If $\mathbb{Z}[X]/I$ is Noetherian, then $\mathbb{F}_2[X]/I_{\mathbb{F}_2}$ is also
Noetherian and by Proposition \ref{bool-f2} it is also a Boolean ring.
From \cite[2.18]{kleiman2013term} it has Krull dimension $0$
and then, it is Artinian, see \cite[Th. 8.5, p.90]{MR0242802}.
Given that $\mathbb{F}_2[X]/I_{\mathbb{F}_2}$ is reduced,
\cite[3.23]{kleiman2013term},
we get that it is isomorphic to $\mathbb{F}_2^n$,
see \cite[Th. 8.7, p.90 and Ex. 28, p.35]{MR0242802}.
Clearly, if $X$ is finite, then $\mathbb{Z}[X]/I$ is Noetherian
(the converse
is false as shown in example $M\omega$ above).
Finally, from Proposition \ref{bool-f2}, $Y\cong\mathbb{F}_2^n$
and from Example \ref{powerset},
$\mathbb{Z}^n\cong \mathbb{Z}[Y]/I\cong\mathbb{Z}[X]/I$.

Assume now that $\rk(\mathbb{Z}[X]/I)<\infty$.
Then, $\mathbb{F}_2[X]/I_{\mathbb{F}_2}\cong\mathbb{F}_2^n$
and again $Y\cong\mathbb{F}_2^n$.
Then, $\mathbb{Z}[X]/I\cong\mathbb{Z}^n$.
\end{proof}

\begin{definition}
Let $X$ be a bounded lattice.
If the rank as a $\mathbb{Z}$-module of $M(X)$
is finite, we set $n(X)=\rk(M(X))$,
else $n(X)=\infty$.
The number $n(X)$ is called the \emph{measurability} of $X$.
Notice that the measurability
completely characterizes measures
on a finite lattice $X$, that is,
$\mathbb{Z}[X]/I\cong \mathbb{Z}^{n}$, $n=n(X)$.

Let $A$ be a field and let $\mathcal{X}_A=\Spec(A[X]/I_A)$.
From Proposition \ref{z-n} it follows that, independently of $A$,
\[
n(X)=\begin{cases}
\#\mathcal{X}_A& \text{if }n(X)<\infty.\\
\infty&\text{if not.}
\end{cases}
\]
\end{definition}

The next result can be found in \cite{geis1}.
\begin{proposition}\label{prod}
Let $X$ and $X'$ be two bounded lattices.
Then, we have a natural isomorphism of rings
\[
\mathbb{Z}[X\times X']/I_{X\times X'}\cong
\mathbb{Z}[X]/I_X\times \mathbb{Z}[X']/I_{X'}.
\]
In particular, $n(X\times X')=n(X)+n(X')$.
\end{proposition}
\begin{proof}
Recall that the lattice structure on $X\times X'$ is coordinate-wise,
\[
(x_1,y_1)\wedge (x_2,y_2)=(x_1\wedge x_2,y_1\wedge y_2),\quad
(x_1,y_1)\vee (x_2,y_2)=(x_1\vee x_2,y_1\vee y_2).
\]
The projection $\pi_1:X\times X'\rightarrow X$
is a lattice map, hence it induces a ring map
$\pi_1:\mathbb{Z}[X\times X']/I_{X\times X'}
\rightarrow \mathbb{Z}[X]/I_X$.
Same for $\pi_2$.
Then,
\[
\pi=\pi_1\times \pi_2:
\mathbb{Z}[X\times X']/I_{X\times X'}\rightarrow
\mathbb{Z}[X]/I_X\times \mathbb{Z}[X']/I_{X'}.
\]
It remains to check that $\pi$ is bijective.
Consider the following linear map,
\[
\psi:\mathbb{Z}^{\oplus X}\times\mathbb{Z}^{\oplus X'}\rightarrow
M(X\times X'),\quad ([x],[y])\mapsto \overline{[(x,y)]}.
\]
It is easy to check that
$\psi(J_X,[0_{X'}])=0$ and
$\psi([0_X],J_{X'})=0$.
Then, $\psi$ induces
$\overline{\psi}:M(X)\times M({X'})\rightarrow M(X\times {X'})$,
the inverse of $\pi$.
\end{proof}

\begin{proposition}\label{noet-com}
Let $X$ be a Boolean lattice.
Then, $n(X)<\infty$ if and only if $X$ is finite.
\end{proposition}
\begin{proof}
If $n(X)<\infty$, then $X\cong\mathbb{F}_2^{n(X)}$ which it is finite.
If $X$ is finite, then $n(X)<\infty$, see Proposition \ref{z-n}.
\end{proof}

\begin{remark}
Let us give a method to compute $n(X)$.
Assume that $f:X\rightarrow X'$ is a surjective lattice
map between two bounded lattices.
Then, it is easy to check that $f_{*}:Y\rightarrow Y'$ is also surjective and then,
$n(X)=n(Y)\ge n(Y')=n(X')$, where $Y$ and $Y'$ are the Boolean lattices associated to $X$ and
$X'$ respectively (see Theorem \ref{teo-Y}).
Hence, it is possible to bound (or compute)
the measurability of a finite lattice
by collapsing edges of its Hasse diagram.
The following list, taken from \cite[Fig.3.5,p.248]{MR2868112},
shows the Hasse diagrams
of all lattices with at most six elements,

\begin{center}
\includegraphics[width=11.7cm]{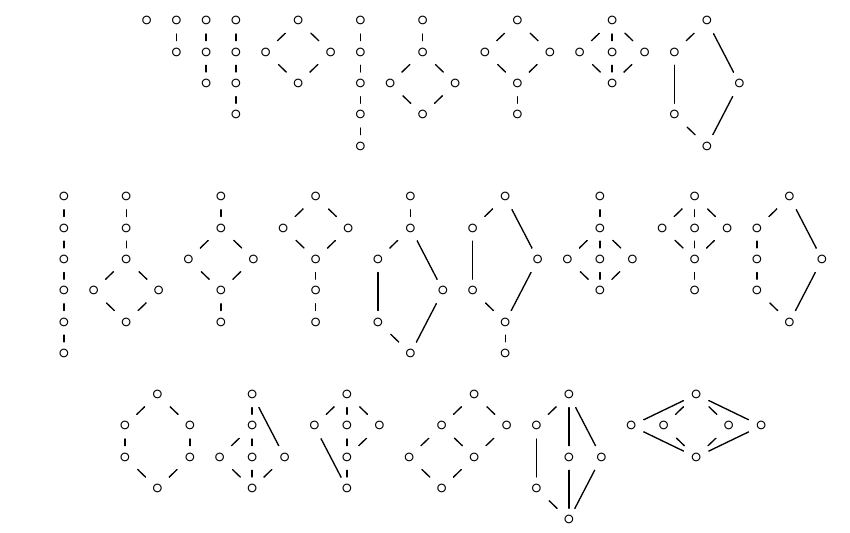}
\end{center}
We computed the corresponding numbers $n(X)$,
\[
0,1,2,3,2,4,3,3,0,2.
\]
\[
5,4,4,4,3,3,1,1,2.
\]
\[
2,1,1,3,0,0.
\]
The $3$ in the last row follows from the fact that the
lattice is the cartesian product of $0\leq 1$ and $0\leq 1\leq 2$.
Then, $3=1+2$.

An effective way to compute the measurability
is by using Gr\"obner bases over $\mathbb{Q}$.
The goal is to find the dimension
of $\mathbb{Q}[X]/I_{\mathbb{Q}}$ as a $\mathbb{Q}$-vector space.
For example, the following lattice has $n(X)=2$,

\begin{center}
\includegraphics[width=3.5cm]{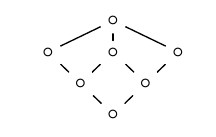}
\end{center}
\end{remark}

\end{document}